\documentclass[11pt,a4paper]{article}
\usepackage[utf8]{inputenc}
\usepackage{amsmath,amssymb,amsfonts,amsthm}
\usepackage{physics}
\usepackage{graphicx}
\usepackage{hyperref}
\usepackage{xcolor}
\usepackage{geometry}
\usepackage{booktabs}
\numberwithin{equation}{section}
\usepackage{lipsum}

\theoremstyle{plain}

\theoremstyle{remark}

\title{\textbf{Bifurcation of Tetrahedral Non-Zonal Flows in the 2D Euler Equations on a Rotating Sphere}}
\author{Yuri Cacchiò\\
\scriptsize{Faculty of Mathematics, University of Vienna, Oskar-Morgenstern-Platz 1, 1090 Vienna, Austria.}\\ \href{mailto:yuri.cacchio@univie.ac.at}{\scriptsize{yuri.cacchio@univie.ac.at}}}
\date{\today}

\begin{document}

\maketitle

\begin{abstract}
We investigate the emergence of finite-amplitude non-zonal flows on the sphere $\mathbb{S}^2$ arising from stationary solutions to the 2D Euler equations. By restricting the Laplace-Beltrami eigenspace to the invariant subspace of the tetrahedral symmetry group $\mathbf{T}$, we bypass the $(2l+1)$-dimensional kernel degeneracy, obtaining a scalar Liapunov-Schmidt reduction. We analyze four distinct physical non-linearities: a polynomial model, the sine-Gordon and sinh-Gordon models, and the exponential (Liouville) model. We explicitly derive the bifurcation parameter via spectral projections, proving that the bifurcation topology (subcritical or supercritical) is not a geometric invariant, but is governed by the parity of the nonlinearity and the mass conservation.

\vspace{0.5cm}
\end{abstract}

\tableofcontents

\section{Introduction}

The 2D Euler equations on the sphere $\mathbb{S}^2$ represent a mathematical framework for large-scale geophysical fluid dynamics \cite{constantin2021modelling,holton2013introduction,taylor2016euler,vallis2017atmospheric}. Let the sphere be parameterized by longitude $\varphi \in [0, 2\pi)$ and latitude $\theta \in \left(-\frac{\pi}{2},\frac{\pi}{2}\right)$, and let $x_3 = \sin\theta$ denote the polar axis. 
The sphere is rotating with constant angular velocity $\omega$ about $x_3$. 

Let $\psi(t, \varphi, \theta)$ be the stream function, such that the relative vorticity is given by the Laplace-Beltrami operator on the unit sphere. In terms of longitude $\varphi$ and latitude $\theta$, this operator is defined as
\begin{equation}
    \Delta\psi = \frac{1}{\cos\theta} \frac{\partial}{\partial\theta} \left( \cos\theta \frac{\partial\psi}{\partial\theta} \right) + \frac{1}{\cos^2\theta} \frac{\partial^2\psi}{\partial\varphi^2}.
\end{equation}
The Euler equation on a rotating frame is expressed as \cite{majda2002vorticity},
\begin{equation}
    \partial_t (\Delta\psi) + J(\psi, \Delta\psi + 2\omega x_3) = 0,
\end{equation}
where $2\omega x_3$ is the planetary Coriolis parameter, and the Jacobian $J$ is defined as 
\begin{equation}
    J(f,g) = \frac{1}{\cos\theta} (\partial_\varphi f \partial_\theta g - \partial_\theta f \partial_\varphi g).
\end{equation}
We are interested in steady-state dynamics ($\partial_t = 0$). Then, the stationary Euler equation reduces to,
\begin{equation}
    J(\psi, \Delta\psi + 2\omega x_3) = 0.
\end{equation}
By classical computations (see for example Arnold and Khesin \cite{arnold1998}), the vanishing of the Jacobian implies that the gradient of the absolute vorticity is everywhere parallel to the gradient of the stream function. As discussed by Constantin and Johnson \cite{constantin2017large}, in regions where the gradient of the stream function does not vanish, the absolute vorticity is a functional of the stream function. This ensures the existence of a sufficiently smooth function $F$ such that,
\begin{equation} \label{eq:CG_euler}
    -\Delta\psi + F(\psi) = 2\omega x_3.
\end{equation}
Equation \eqref{eq:CG_euler} admits 1D \textit{zonal flows}, which depend only on the latitudinal coordinate $x_3$. However, at critical parameter thresholds, these symmetric states can exhibit spontaneous symmetry breaking, giving rise to fully 2D, non-zonal \textit{Rossby-Haurwitz waves} \cite{haurwitz1940motion,rossby1939relations} via a local bifurcation. 

To isolate the topological nature of this bifurcation without the geometric asymmetry introduced by the Coriolis force, we follow the idea of Constantin and Germain \cite[Sec. 6.1]{constantin2022} and investigate the non-rotating limit. By setting the planetary rotation to zero, that is $\omega = 0$, the stationary Euler equation \eqref{eq:CG_euler}, simplifies to the following semilinear elliptic problem,
\begin{equation} \label{eq:general}
    -\Delta\psi + F(\lambda,\psi) = 0 \quad \text{on } \mathbb{S}^2,
\end{equation}
where $F(\lambda, 0) = 0$, and $\lambda \in \mathbb{R}$ acts as the bifurcation parameter governing the nonlinearity. 

Returning to the rotating sphere ($\omega > 0$), the planetary rotation breaks the global symmetry. Consequently, constructing exact solutions becomes challenging. For instance, while Stuart-type vortices can be formulated explicitly on a fixed sphere \cite{crowdy2004}, modeling them on a rotating sphere typically requires perturbative approaches around background zonal flows \cite{constantin2021stuart}.
Despite this broken symmetry, a correspondence exists between stationary solutions at $\omega = 0$ and traveling waves at $\omega > 0$. As shown by Constantin and Germain \cite{constantin2022}, any stationary solution $\psi_0(\varphi, \theta)$ of \eqref{eq:general} maps to a solution of the rotating Euler equations via the transformation $\psi_\omega(\varphi, \theta, t) = \psi_0(\varphi + \omega t, \theta) + \omega x_3$. 
This mapping transforms the stationary states into non-zonal waves that propagate westward at a phase speed $c = \omega$ relative to the planetary surface.

Recent studies on the 2D Euler equations on $\mathbb{S}^2$ have addressed the bifurcation of vortex caps from zonal flows \cite{garcia2025dynamics}, the regularization of point vortices via Lyapunov-Schmidt reduction \cite{sakajo2025regularization}, and the rigidity of stratospheric travelling waves \cite{constantin2025rigidity}.

Our objective is to determine how the profile function $F$ characterizes the dynamics of the bifurcating tetrahedral non-zonal flows. In particular, we show that the bifurcation topology (subcritical or supercritical) is highly sensitive to the nonlinear interactions generated by the specific nonlinearity.

To rigorously investigate this bifurcation, one must be careful of the geometric degeneracy of the spherical Laplacian (see Sec \ref{sec:dim_prob} below). The multidimensionality of the eigenspaces on the sphere prohibits the direct application of the classical one-dimensional bifurcation theorems. To bypass this problem, we follow \cite{constantin2022} by searching for solutions within a restricted geometric class, the isotropy subgroup of tetrahedral symmetry (see Sec \ref{Sec:Tetra_Group}). 

The tetrahedral group provides a strategic framework \cite{golubitsky1988singularities}. In fact, among the polyhedral symmetries, it is the lowest-degree configuration and corresponds to spherical harmonics of degree $l=3$. This restriction removes the geometric degeneracy, collapsing the problem into a one-dimensional kernel.

While the existence of these bifurcating branches can be established using the Rabinowitz global bifurcation approach \cite{rabinowitz1971some}(as proved by Constantin and Germain \cite{constantin2022}), such theorems do not take into account the local geometry and direction of the branch. More in detail, to understand whether the non-zonal flows emerge subcritically or supercritically requires a precise asymptotic expansion. 

In this paper, we show that the nature of these bifurcations is sensitive to the structure of the specific nonlinearity $F(\lambda, \psi)$. By applying the Liapunov-Schmidt reduction \cite{kielhofer2012bifurcation} to four distinct physically relevant models, a generalized polynomial \cite{pedlosky1987geophysical}, the sine-Gordon and sinh-Gordon equations \cite{flierl1987isolated,stuart1967finite}, and the Liouville exponential model \cite{joyce1973negative}, we understand how the convexity, parity, and geometric constraints of each function affect the nonlinear interactions. Finally, we prove that this nonlinear interaction between the non-zonal flow and the even-degree harmonics governs the topology of the bifurcation.

The organization of this paper is as follows. In Section \ref{sec:Liap-S}, we formulate the general perturbative framework using the Liapunov-Schmidt reduction and highlight the geometric degeneracy of the spherical kernel. Section \ref{Sec:Tetra_Group} deals with the restriction of the functional space to the invariant subspace of the tetrahedral symmetry group $\mathbf{T}$ to isolate a one-dimensional kernel. In Section \ref{sec:universal_order1}, we characterize the universal leading-order structure of the bifurcating non-zonal flow ($l=3$). In Section \ref{sec:bifurcation}, we investigate the nonlinear interactions for four distinct profile functions: generalized polynomials (Section \ref{sec_polyn}), the sine-Gordon model (Section \ref{sec:sine}), the sinh-Gordon model (Section \ref{sec:sinh}), and the Liouville exponential model (Section \ref{sec:exp}). Finally, Section \ref{sec:conclusions} synthesizes our results and provides concluding remarks.

\section{The Liapunov-Schmidt Reduction}\label{sec:Liap-S}

Before addressing specific nonlinearities, we formulate a general perturbative framework. Let us define the nonlinear operator $\mathcal{G}: \mathbb{R} \times H^s(\mathbb{S}^2) \to H^{s-2}(\mathbb{S}^2)$ as,
\begin{equation}\label{eq:G}
    \mathcal{G}(\lambda, \psi) := -\Delta\psi + F(\lambda, \psi).
\end{equation}
Without loss of generality, we assume that $F(\lambda, 0) = 0$ for all $\lambda$, such that $\mathcal{G}(\lambda, 0) = 0$ defines a trivial solution branch. 
A necessary condition for spontaneous symmetry breaking (bifurcation) to occur at a critical parameter $\lambda^*$ along the trivial branch is the failure of the Implicit Function Theorem. This mathematically requires that the linearized operator around the state $(\lambda^*, 0)$ has a non-trivial kernel.
The linearized operator around $(\lambda^*, 0)$ is defined as the Fréchet derivative \begin{equation}
    L := \partial_\psi \mathcal{G}(\lambda^*, 0).
\end{equation}
In particular, taking the variational derivative of \eqref{eq:G} yields,
\begin{equation}
    L = -\Delta + \partial_\psi F(\lambda^*, 0) I.
\end{equation}
In order to track the non-trivial solutions in the neighborhood of this bifurcation point, we use the Liapunov-Schmidt reduction. Since $L$ is an elliptic, self-adjoint linear differential operator on a compact manifold ($\mathbb{S}^2$), classical elliptic theory \cite{kielhofer2012bifurcation} tells us that $L$ is a Fredholm operator of index zero with a closed range. 

Notice that the self-adjointness property, $L = L^*$, allows us to apply the Fredholm Alternative. 
This results in the orthogonal decomposition of our functional space to the direct sum of the kernel and the image of the operator,
\begin{equation}
    H^s(\mathbb{S}^2) = \text{Ker}(L) \oplus \text{Im}(L),
\end{equation}
with,
\begin{equation}\label{eq:lyap_cond}
    \text{Im}(L) = (\text{Ker}(L))^\perp.
\end{equation}
Let $P$ be the $L^2$-orthogonal projection operator onto $\ker(L)$, and let $Q = (I - P)$ be the projection in $\text{Im}(L)$. Any function $\psi$ can be uniquely decomposed as,
\begin{equation}
    \psi = v + w \quad \text{where} \quad \begin{cases}
           v \in \ker(L), \\ w \in \text{Im}(L).
    \end{cases}
\end{equation}
Applying the projectors $P$ and $Q$ to the equation $\mathcal{G}(\lambda, v+w) = 0$, we derive an equivalent coupled system,
\begin{align}
    Q \mathcal{G}(\lambda, v+w) &= 0, \quad \text{(Auxiliary Equation)} \label{eq:aux}\\
    P \mathcal{G}(\lambda, v+w) &= 0. \quad \text{(Bifurcation Equation)} \label{eq:bifurc}
\end{align}
Let us emphasize that $Q L = L|_{\text{Im}(L)}$ is, by construction, an isomorphism from $\text{Im}(L)$ to itself. Therefore, by the Implicit Function Theorem, the Auxiliary Equation \eqref{eq:aux} can be uniquely solved for $w$ in a sufficiently small neighborhood of $(\lambda^*, 0)$. Moreover, the solution $w = w(\lambda, v)$ constructed in this way satisfies $w(\lambda^*, 0) = 0$.

Substituting this map back into Equation \eqref{eq:bifurc}, the original problem is reduced to a finite-dimensional problem on the kernel,
\begin{equation}\label{eq:bif}
    \Phi(\lambda, v) := P \mathcal{G}(\lambda, v + w(\lambda, v)) = 0.
\end{equation}
To extract the local structure of the bifurcating branch, we parameterize the solutions using a small parameter $\varepsilon$. Setting $v = \varepsilon \psi_1$ (where $\psi_1$ spans $\ker(L)$), we expand the fields asymptotically:
\begin{align}
    \psi &= \varepsilon \psi_1 + \varepsilon^2 \psi_2 + \varepsilon^3 \psi_3 + \mathcal{O}(\varepsilon^4),\label{eq:psi} \\
    \lambda &= \lambda^* + \varepsilon \lambda_1 + \varepsilon^2 \lambda_2 + \mathcal{O}(\varepsilon^3),\label{eq:lambda}
\end{align}
where $\psi_k \in \text{Im}(L)$ for $k \ge 2$. Substituting these series in $\mathcal{G}(\lambda, \psi) = 0$ transforms the original equation into a sequence of linear problems, solvable order by order.

\subsection{The Multi-Dimensional Kernel}\label{sec:dim_prob}
The solution of the reduced function $\Phi(\lambda, v) = 0$ depends on the dimension of $\ker(L)$. If $\dim(\ker(L)) = 1$, the function $\Phi$ is a scalar equation, and the classical Crandall-Rabinowitz theorem \cite{kielhofer2012bifurcation} guarantees a unique bifurcating branch. 

However, on the sphere $\mathbb{S}^2$ the eigenvalues of the Laplace-Beltrami operator corresponding to a spherical harmonic degree $l$ are non-simple. In particular,
\begin{equation}
    \dim(\ker(L)) = 2l + 1.
\end{equation}
For $l \ge 1$, the Bifurcation Equation \eqref{eq:bif} becomes a system of $2l+1$ nonlinear equations. For example, bifurcating at $l=3$ gives a problem with $7$ equations. Extracting information directly from this system without any other mathematical structure becomes challenging.

\section{The Tetrahedral Group}\label{Sec:Tetra_Group}
To bypass the degeneracy of the $(2l+1)$ dimensional kernel, we exploit the symmetries of the system. The domain $\mathbb{S}^2$, the Laplace-Beltrami operator $\Delta$, and the local nonlinearities $F(\lambda, \psi)$ considered in this work are all invariant under the action of the orthogonal group $O(3)$ (the group of all distance-preserving transformations of $\mathbb{R}^3$). Therefore, the nonlinear operator $\mathcal{G}$ commutes with $O(3)$, implying that the symmetries of the solutions are preserved.

We restrict our functional space to a specific isotropy subgroup. Let $\mathbf{T} \subset O(3)$ denote the tetrahedral group, the symmetry group of a regular tetrahedron inscribed in the sphere. We define the invariant subspace $X_{\mathbf{T}}$ as,
\begin{equation}
    X_{\mathbf{T}} = \{ \psi \in H^s(\mathbb{S}^2) \mid \psi(g \cdot x) = \psi(x) \quad \forall g \in \mathbf{T}, \forall x \in \mathbb{S}^2 \}.
\end{equation}
By the Principle of Symmetric Criticality \cite{palais1979} and standard equivariant bifurcation lemmas \cite{golubitsky1988singularities}, solutions in the restricted subspace $X_{\mathbf{T}}$ are mathematically well-posed. Since the operator $\mathcal{G}$ is equivariant, it  maps $X_{\mathbf{T}}$ to itself $\mathcal{G}: X_{\mathbf{T}} \to X_{\mathbf{T}}$. This ensures that any local solution found within this symmetric subspace is an exact solution of the full Euler equations \eqref{eq:general}.
Indeed, the above symmetry forces all orthogonal projections onto the remaining components of the $(2l+1)$-dimensional kernel to vanish.

The advantage of working in the restricted subspace $X_{\mathbf{T}}$ is related to the group representation theory \cite{golubitsky1988singularities}. Under the action of the tetrahedral group $\mathbf{T}$, the 7-dimensional representation of $O(3)$ for the $l=3$ eigenspace decomposes. This decomposition contains exactly \textit{one} fully invariant representation. 

Consequently, within $X_{\mathbf{T}}$, the previously 7-dimensional kernel collapses to dimension one,
\begin{equation}
    \dim(\ker(L|_{X_{\mathbf{T}}})) = 1.
\end{equation}
This reduction to a single Bifurcation Equation satisfies the requisite of the local Crandall-Rabinowitz theorem.

The unique generator is the real tetrahedral harmonic \cite{muller1966spherical},
\begin{equation} \label{eq:Ystar}
    Y^*(\varphi, \theta) = R_3^{-2}(\varphi, \theta) = \sqrt{\frac{105}{16\pi}} \sin\theta \cos^2\theta \cos(2\varphi).
\end{equation}
Then, the restricted kernel is defined as,
\begin{equation}
    \ker(L|_{X_{\mathbf{T}}}) = \text{span}\{Y^*\}.
\end{equation}
Having this one-dimensional kernel, allow us to use the Liapunov-Schmidt framework introduced in Section \ref{sec:Liap-S}.

\section{Universal Leading-Order}\label{sec:universal_order1}

Before specifying the form of the nonlinearities, we rigorously determine the leading-order structure of the bifurcating solutions. For any generic smooth nonlinearity satisfying the trivial-branch condition $F(\lambda, 0) = 0$, the Taylor expansion with respect to $\psi$ evaluated at the bifurcation point $(\lambda^*, 0)$ is given by,
\begin{equation}
    F(\lambda, \psi) = \partial_\psi F(\lambda^*, 0) \psi + \mathcal{O}(\psi^2).
\end{equation}
Substituting the perturbative series 
\begin{align*}
    \psi &= \varepsilon \psi_1 + \mathcal{O}(\varepsilon^2),\\
    \lambda &= \lambda^* + \mathcal{O}(\varepsilon),
\end{align*}
into equation \eqref{eq:general}, and isolating the terms of order $\mathcal{O}(\varepsilon^1)$, the problem reduces to the linearized equation,
\begin{equation}\label{eq:oreder1}
    -\Delta\psi_1 + \partial_\psi F(\lambda^*, 0) \psi_1 = 0,
\end{equation}
that is,
\begin{equation}
    L\psi_1 = 0,
\end{equation}
where $L$ is the Fréchet derivative defined in Section \ref{sec:Liap-S}.
This implies that $$\psi_1\in \ker(L).$$ Since we work in the invariant subspace $X_{\mathbf{T}}$, the kernel is one-dimensional and spanned by the tetrahedral harmonic \eqref{eq:Ystar}. By absorbing the constant into the parameter $\varepsilon$, we deduce,
\begin{equation}
    \psi_1 = Y^*.
\end{equation}
Then, to leading order, the bifurcating stream function exhibits a universal geometry, regardless of the specific structure of the higher-order nonlinearities, i.e.,

 \begin{equation}\label{eq:universal_psi}
    \psi(\varphi, \theta) \sim \pm \varepsilon \left( \sqrt{\frac{105}{16\pi}} \sin\theta \cos^2\theta \cos(2\varphi) \right) + \mathcal{O}(\varepsilon^2).
\end{equation}
From this explicit formulation, we deduce the following physical properties of the non-zonal flow.
\begin{itemize}
    \item \textbf{Mid-latitudes Dynamic:} The velocity field cannot have any vortex at the equator, $\theta=0$, and at the poles, $\theta=\pm\pi/2$. In fact, $\sin(0)=0$ and $\cos(\pm\pi/2)=0$ force \eqref{eq:universal_psi} to vanish at these points. Therefore, the fluid action is entirely confined to mid-latitudes. 
    \item \textbf{Tetrahedral Vortex:} The stream function exhibits four positive maxima and four negative minima. Optimizing the latitudinal profile $f(\theta) = \sin\theta\cos^2\theta$ by setting the derivative $f'(\theta)=0$, we obtain,
    \begin{equation}
        \cos^3\theta - 2\sin^2\theta\cos\theta = 0,
    \end{equation}
    from which,
    \begin{equation}
         \tan\theta = \pm \frac{1}{\sqrt{2}}.
    \end{equation}
    Thus, the critical latitudes for the vortex are fixed at $\theta \approx \pm 0.615$ rad (approximately $\pm 35.26^\circ$). Combined with $\cos(2\varphi)$ at $0^\circ, 90^\circ, 180^\circ,$ and $270^\circ$, the positive vortices map the vertices of a regular tetrahedron inscribed in the sphere, while the negative vortices perfectly map the dual tetrahedron.
    \item \textbf{Equatorial Asymmetry:} Due to the $\sin\theta$ term, \eqref{eq:universal_psi} is odd with respect to the equator,
    \begin{equation}
        Y^*(-\theta) = -Y^*(\theta).
    \end{equation}
     At leading order, the vortices exhibit opposite signs and intensities between the Northern and Southern hemispheres.
\end{itemize}
As we demonstrate, while this leading-order spatial structure is universal at $\mathcal{O}(\varepsilon)$, the second order $\mathcal{O}(\varepsilon^2)$ depends strictly on the fixed model. Indeed, the nonlinearity determines whether this hemispheric anti-symmetry is preserved or broken by the nonlinear interactions.

\section{Bifurcation \label{sec:bifurcation}}
The compact geometry of the domain imposes a constraint on the nonlinear flow. Integrating equation \eqref{eq:general} over the sphere $\mathbb{S}^2$, the surface integral of the Laplacian vanishes identically by Gauss's Divergence Theorem,
\begin{equation}
    \iint_{\mathbb{S}^2} \Delta\psi \, d\sigma = 0.
\end{equation}
Consequently, any admissible function $F$ must satisfy the geometric constraint,
\begin{equation} \label{eq:gauss_constraint}
    \iint_{\mathbb{S}^2} F(\lambda, \psi) \, d\sigma = 0.
\end{equation}
Notice that the mechanism required to satisfy this constraint strictly depends on the parity of the nonlinearity. In this section, we deal with four different physical models: a polynomial model \cite{pedlosky1987geophysical}, the sine-Gordon and sinh-Gordon models \cite{flierl1987isolated, stuart1967finite}, and the exponential (Liouville) model \cite{joyce1973negative}. 
Before proceeding with the bifurcation analysis, we briefly outline how the four chosen models satisfy condition \eqref{eq:gauss_constraint}.
\begin{itemize}
    \item \textbf{Polynomial Nonlinearities (Section \ref{sec_polyn}):} For generalized polynomials, the quadratic term $\psi^2$ generates a strictly positive mass. To prevent the violation of \eqref{eq:gauss_constraint}, the system must generate a non-zero spatial mean (the $l=0$ zonal mode) in the higher-order corrections to absorb this integral.
    
    \item \textbf{Odd Nonlinearities (Sections \ref{sec:sine} and \ref{sec:sinh}):} 
    For odd functions $F(-\psi) = -F(\psi)$, such as the sine-Gordon and sinh-Gordon equations, the constraint is trivially satisfied. The anti-symmetry of the solutions ensures \eqref{eq:gauss_constraint}, forbidding the generation of zonal flows.
    
    \item \textbf{Positive Functions (Section \ref{sec:exp}):} For positive nonlinearities like the Liouville exponential model, the Gauss constraint \eqref{eq:gauss_constraint} is not satisfied. This forces the introduction of an integral term (a topological mass shift) to match condition \eqref{eq:gauss_constraint}.
\end{itemize}
The geometric framework is now fully established. As demonstrated in Section \ref{sec:universal_order1}, the primary spatial kernel is fixed to the tetrahedral harmonic $\psi_1 = Y^*$, corresponding to the Laplacian eigenvalue $l(l+1)=12$. We now perform the local asymptotic expansion up to $\mathcal{O}(\varepsilon^3)$ for each specific regime to determine their critical parameters $\lambda^*$ and bifurcation topologies.

\subsection{The Polynomial Model}\label{sec_polyn}
Following Theorem 9 and Remark 2 in \cite{constantin2022}, we first consider a polynomial model. We fix
\begin{equation}
    P(\lambda) = \mu_1\lambda^3 - [\mu+12]\lambda
\end{equation}
where $\mu, \mu_1 > 0$. Then, the nonlinearity becomes
\begin{align}\label{eq:pol_F}
    F(\lambda,\psi) &= P(\lambda+\psi) - P(\lambda) \nonumber \\
    &=\mu_1\left(3\lambda\psi^2 + \psi^3\right) + \left(3\mu_1\lambda^2 - [\mu+12]\right)\psi.
\end{align}
We now substitute the asymptotic expansions \eqref{eq:psi} and \eqref{eq:lambda} into the general equation \eqref{eq:general}, with $F(\lambda,\psi)$ defined as in \eqref{eq:pol_F}. By collecting terms by powers of $\varepsilon$, we get 
\subsubsection*{Order $\mathcal{O}(\varepsilon^1)$}
The universal condition derived in Section \ref{sec:universal_order1} $$\partial_\psi F(\lambda^*, 0) \psi_1 = -12\psi_1,$$ translates to $$3\mu_1(\lambda^*)^2 - \mu - 12 = -12.$$ Solving for the unknown parameter  $\lambda^*$, we find the bifurcation point,
\begin{equation}\label{eq:bifurc_cond}
    \lambda^* = \sqrt{\frac{\mu}{3\mu_1}}.
\end{equation}

\subsubsection*{Order $\mathcal{O}(\varepsilon^2)$}
Collecting terms of order $\mathcal{O}(\varepsilon^2)$ and recalling the notation $ L = -\Delta + \partial_\psi F(\lambda^*, 0) I$, we derive
\begin{equation}\label{eq:second_order_pol}
    L\psi_2 + 6\mu_1\lambda^*\lambda_1\psi_1 + 3\mu_1\lambda^*\psi_1^2 = 0.
\end{equation}
We now project equation \eqref{eq:second_order_pol} onto the kernel spanned by $\psi_1 = Y^*$.
Since $L$ is a self-adjoint operator, the first term vanishes
$\langle L\psi_2, Y^* \rangle = \langle \psi_2, LY^* \rangle = 0$, resulting in 
\begin{equation}\label{eq:pol_O2}
    6\mu_1\lambda^*\lambda_1 \|Y^*\|_{L^2}^2 + 3\mu_1\lambda^* \iint_{\mathbb{S}^2} (Y^*)^3 \, d\sigma = 0.
\end{equation}
Notice that $Y^*$ is an odd-degree spherical harmonic ($l=3$), then $(Y^*)^3$ remains an odd function over the sphere. It follows that the integral term in \eqref{eq:pol_O2} vanishes. This property forces to
\begin{equation}
    \lambda_1 = 0.
\end{equation}
While the parameter correction $\lambda_1$ vanishes, the spatial correction $\psi_2$ is non-zero. Therefore, equation \eqref{eq:second_order_pol}  reduces to,
\begin{equation}
    L\psi_2 =- 3\mu_1\lambda^*(Y^*)^2,
\end{equation}
from which we obtain,
\begin{equation}\label{eq:psi_2}
    \psi_2 = -3\mu_1\lambda^* L^{-1}((Y^*)^2).
\end{equation}
Since the source term $(Y^*)^2$ is non-negative and of even parity, its spectral decomposition is spanned by even spherical harmonics, $l \in \{0, 2, 4, 6\}$. As discussed above, the $l=0$ projection is related to the Gauss constraint \eqref{eq:gauss_constraint} to balance the positive mass of the quadratic term.
Furthermore, because these even harmonics are orthogonal to the odd $l=3$ kernel, the application of the inverse operator $L^{-1}$ is well-posed.

\subsubsection*{Order $\mathcal{O}(\varepsilon^3)$}
Recalling that $\lambda_1 = 0$, at the third order the expansion gives,
\begin{equation}
    L\psi_3 + 6\mu_1\lambda^*\lambda_2\psi_1 + 6\mu_1\lambda^*\psi_1\psi_2 + \mu_1\psi_1^3 = 0.
\end{equation}
Let us project this equation in the kernel spanned by $Y^*$. We note that $L\psi_3$ vanishes since $L$ is self-adjoint (as $\psi_2$ above). Since $\|Y^*\|^2 = 1$, we  isolate the second-order bifurcation parameter
\begin{equation} \label{eq:lambda2_raw}
    \lambda_2 = -\frac{1}{6\lambda^*} \left[ \int_{\mathbb{S}^2} (Y^*)^4 \, d\sigma + 6\lambda^* \int_{\mathbb{S}^2} (Y^*)^2 \psi_2 \, d\sigma \right].
\end{equation}

\subsubsection*{Sign of $\lambda_2$}
We evaluate the non-local integral involving $\psi_2$. Substituting \eqref{eq:psi_2} into \eqref{eq:lambda2_raw}, we obtain,
\begin{align} \label{eq:resolvent_substitution}
    6\lambda^* \int_{\mathbb{S}^2} (Y^*)^2 \psi_2 \, d\sigma &= 6\lambda^* \int_{\mathbb{S}^2} (Y^*)^2 \left[ -3\mu_1\lambda^* L^{-1}((Y^*)^2) \right] \, d\sigma \nonumber \\
    &= -18\mu_1(\lambda^*)^2 \int_{\mathbb{S}^2} (Y^*)^2 L^{-1}((Y^*)^2) \, d\sigma.
\end{align}
We rewrite the coefficient $18\mu_1(\lambda^*)^2$ by using the linear bifurcation condition \eqref{eq:bifurc_cond} derived at order $\mathcal{O}(\varepsilon^1)$, from which 
\begin{equation}
    18\mu_1 \left( \frac{\mu}{3\mu_1} \right) = 6\mu.
\end{equation}
Thus, equation \eqref{eq:lambda2_raw} becomes,
\begin{equation} \label{eq:lambda2_reduced}
    \lambda_2 = -\frac{1}{6\lambda^*} \left[ \int_{\mathbb{S}^2} (Y^*)^4 \, d\sigma - 6\mu \int_{\mathbb{S}^2} (Y^*)^2 L^{-1}((Y^*)^2) \, d\sigma \right].
\end{equation}
To evaluate these integrals, we decompose the term $(Y^*)^2$ into the real orthonormal spherical harmonic basis $\{R_l^m\}$ \cite{muller1966spherical},
\begin{equation}
    (Y^*)^2 = \sum_{l \in \{0, 2, 4, 6\}} \sum_{m \in \{0, 4\}} c_{l,m} R_l^m.
\end{equation}
Let us evaluate the projection onto the zonal mode $l=2$. The corresponding normalized harmonic is $R_2^0(\theta) \sim (3\sin^2\theta - 1)$. The projection integral is proportional to
\begin{equation}
    c_{2,0} \sim \int_{-\pi/2}^{\pi/2} \sin^2\theta \cos^5\theta (3\sin^2\theta - 1) \, d\theta.
\end{equation}
Using the change of variable $x = \sin\theta$ ($dx = \cos\theta \, d\theta$) and $\cos^4\theta = (1-x^2)^2$, the integral becomes the following polynomial
\begin{align*}
    \int_{-1}^1 x^2 (1-x^2)^2 (3x^2 - 1) \, dx &= \int_{-1}^1 (3x^8 - 7x^6 + 5x^4 - x^2) \, dx \\
    &= 2 \left( \frac{3}{9} - \frac{7}{7} + \frac{5}{5} - \frac{1}{3} \right) = 0.
\end{align*}
Thus, 
\begin{equation}\label{eq:c20}
    c_{2,0} = 0.
\end{equation}
We apply Parseval's identity to the first term in \eqref{eq:lambda2_reduced},
\begin{equation}\label{eq:parseval}
    \int_{\mathbb{S}^2} (Y^*)^4 \, d\sigma = \sum c_{l,m}^2.
\end{equation}
Thereafter, we use the spectral representation of the inverse operator $L^{-1} = (-\Delta - 12)^{-1}$, whose eigenvalues are $(l(l+1)-12)^{-1}$
so that we write $\lambda_2$ as
\begin{align}
    \lambda_2 &= -\frac{1}{6\lambda^*} \sum_{l \in \{0, 4, 6\}} \sum_{m \in \{0, 4\}} c_{l,m}^2 \left( 1 - \frac{6\mu}{l(l+1) - 12} \right)\notag\\
     &= -\frac{1}{6\lambda^*} \sum_{l \in \{0, 4, 6\}} \sum_{m \in \{0, 4\}} c_{l,m}^2 \left( 1 + \frac{6\mu}{12 - l(l+1)} \right).\label{eq:lambda_2_ris}
\end{align}
Let us emphasize that the low-frequency zonal mode $l=0$ yields a strictly positive term ($12-0=12$). On the other hand, the high-frequency modes $l=4, 6$ yield negative denominators ($12-20=-8$ and $12-42=-30$).

To rigorously evaluate the overall sign, we use the explicit spectral weights,
\begin{equation*}
    c_{0,0}^2 = \frac{1}{4\pi}, \quad \sum_{m \in \{0, 4\}} c_{4,m}^2 = \frac{189}{1936\pi}, \quad \sum_{m \in \{0, 4\}} c_{6,m}^2 = \frac{5111}{25168\pi}.
\end{equation*}
Evaluating the sum in \eqref{eq:lambda_2_ris}, we get
\begin{equation}\label{eq:nonlin_pol}
    \sum_{l,m} \frac{c_{l,m}^2}{12 - l(l+1)} = \frac{1}{12}\left(\frac{1}{4\pi}\right) - \frac{1}{8}\left(\frac{189}{1936\pi}\right) - \frac{1}{30}\left(\frac{5111}{25168\pi}\right) = +\frac{5621}{3020160\pi} > 0.
\end{equation}
Since both \eqref{eq:parseval} and \eqref{eq:nonlin_pol} are strictly positive quantities and we have a minus sign in \eqref{eq:lambda_2_ris}, we conclude that,
\begin{equation}
    \lambda_2 < 0.
\end{equation}
This proves that the bifurcation is \textit{subcritical}, providing non-trivial solutions for $\lambda < \lambda^*$. The local asymptotic $\psi$ scales as,
\begin{equation}
    \psi \sim \pm \sqrt{\frac{\lambda^* - \lambda}{|\lambda_2|}} \left( \sqrt{\frac{105}{16\pi}} \sin\theta \cos^2\theta \cos(2\varphi) \right) + \mathcal{O}(\lambda^* - \lambda).
\end{equation}

\subsection{The Sine-Gordon Model}\label{sec:sine}
Let us now consider the sine-Gordon model on the sphere,
\begin{equation}
    -\Delta\psi - \lambda \sin(\psi) = 0.
\end{equation}
By Taylor expanding the nonlinearity, we obtain an asymptotic expansion exclusively in terms of odd powers,
\begin{equation}\label{eq:taylor_sin}
    -\lambda \sin(\psi) = -\lambda \left( \psi - \frac{1}{6}\psi^3 + \mathcal{O}(\psi^5) \right).
\end{equation}
We now apply the same perturbative procedure as for the polynomial case (Sec. \ref{sec_polyn}), setting $F(\lambda, \psi) = -\lambda\sin(\psi)$.

\subsubsection*{Order $\mathcal{O}(\varepsilon^1)$}
As established in Section \ref{sec:universal_order1}, the first order is given by \eqref{eq:oreder1}. For the sine-Gordon model, we compute 
\begin{equation}
    \partial_\psi (\lambda \sin\psi)|_{(\lambda=\lambda^*,\psi=0)} = \lambda^* \cos(0) = \lambda^*.
\end{equation}
Then, the linearized operator reduces to the shifted Laplacian, 
\begin{equation}
    L = -\Delta - \lambda^* I.
\end{equation}
Recalling that $\psi_1 = Y^*$ is an eigenfunction of the Laplacian with eigenvalue $-l(l+1)=-12$, the condition $L\psi_1 = 0$ yields
\begin{equation}
    \lambda^* = 12.
\end{equation}

\subsubsection*{Order $\mathcal{O}(\varepsilon^2)$}
Due to the absence of quadratic terms $\psi^2$ in the Taylor expansion \eqref{eq:taylor_sin}, the second-order equation simplifies to,
\begin{equation}\label{eq:sine_order2}
    L\psi_2 - \lambda_1 Y^* = 0.
\end{equation}
We project the equation onto the kernel spanned by $Y^*$. By the self-adjoint property of $L$, the term $\langle L\psi_2, Y^* \rangle = \langle \psi_2, L Y^* \rangle$ vanishes identically, resulting in 
\begin{equation}
    -\lambda_1 \|Y^*\|_{L^2}^2 = 0,
\end{equation}
from which,
\begin{equation}
    \lambda_1 = 0.
\end{equation}
Equation \eqref{eq:sine_order2} therefore reduces to 
\begin{equation}
    L\psi_2 = 0.
\end{equation}
Imposing the standard Lyapunov-Schmidt orthogonality condition $\psi_2 \perp \ker(L)$, the only admissible solution is trivial,
\begin{equation}
    \psi_2 = 0.
\end{equation}
The odd parity of the vector field prevents the formation of the $l=0$ zonal flow.

\subsubsection*{Order $\mathcal{O}(\varepsilon^3)$} 
Recalling that $\lambda_1 = 0$ and $\psi_2 = 0$, at the third order the expansion gives,
\begin{equation}
    L\psi_3 - \lambda_2 \psi_1 + 12 \left( \frac{\psi_1^3}{6} \right) = 0.
\end{equation}
Once again, projecting this equation onto the primary wave $Y^*$ and exploiting the self-adjointness of $L$, we derive the bifurcation equation,
\begin{equation}
    -\lambda_2 \|Y^*\|^2 + 2 \int_{\mathbb{S}^2} (Y^*)^4 d\sigma = 0,
\end{equation}
from which,
\begin{equation}
     \lambda_2 = 2 \int_{\mathbb{S}^2} (Y^*)^4 d\sigma.
\end{equation}
Since the integral of a positive function is strictly positive, we conclude that 
\begin{equation}
    \lambda_2 > 0.
\end{equation}
The bifurcation is \textit{supercritical}. Solutions exist for $\lambda > 12$,
\begin{equation}
    \psi \sim \pm \sqrt{\frac{\lambda - \lambda^*}{\lambda_2}} \left( \sqrt{\frac{105}{16\pi}} \sin\theta \cos^2\theta \cos(2\varphi) \right) + \mathcal{O}(\lambda - \lambda^*).
\end{equation}

\subsection{The sinh-Gordon Model}\label{sec:sinh}
In this section, we investigate the sinh-Gordon model on the sphere,
\begin{equation}
    -\Delta\psi - \lambda \sinh(\psi) = 0.
\end{equation}
As for the sine-Gordon equation, this nonlinearity exhibits odd parity. However, its local Taylor expansion changes the sign of the leading-order nonlinear term, yielding a positive cubic coefficient,
\begin{equation}
    \sinh(\psi) = \psi + \frac{1}{6}\psi^3 + \mathcal{O}(\psi^5).
\end{equation}
As above, we now separate the perturbative orders.
\subsubsection*{Order $\mathcal{O}(\varepsilon^1)$} 
Similarly to the sin-Gordon model (Sec. \ref{sec:sine}), noting that
\begin{equation}
    \partial_\psi (\lambda \sinh\psi)|_{(\lambda=\lambda^*,\psi=0)} = \lambda^* \cosh(0) = \lambda^*,
\end{equation}
we compute the bifurcation point
\begin{equation}
    \lambda^* = 12.
\end{equation}

\subsubsection*{Order $\mathcal{O}(\varepsilon^2)$} 
Exactly as in the sine-Gordon model, the absence of quadratic terms guarantees that 
$ \lambda_1 = 0$ and forces $\psi_2 \equiv 0$.

\subsubsection*{Order $\mathcal{O}(\varepsilon^3)$} 
Recalling that $ \lambda_1 = 0$ and $\psi_2 \equiv 0$, the equation is given by 
\begin{equation}
      L\psi_3 - \lambda_2 \psi_1 - 2\psi_1^3 = 0.
\end{equation} 
Projecting onto the wave $\psi_1 = Y^*$, we obtain
\begin{equation}
    -\lambda_2 \|Y^*\|_{L^2}^2 - 2 \int_{\mathbb{S}^2} (Y^*)^4 \, d\sigma = 0.
\end{equation}
Since $ \|Y^*\|_{L^2}^2=1$, we isolate the second-order bifurcation parameter,
\begin{equation}
    \lambda_2 = -2 \int_{\mathbb{S}^2} (Y^*)^4 \, d\sigma.
\end{equation}
The integral of a non-trivial fourth power is positive, and the minus sign forces
\begin{equation}
    \lambda_2 < 0.
\end{equation}
We conclude that the bifurcation of the sinh-Gordon model is \textit{subcritical}, with non-trivial solutions for $\lambda < 12$. Unlike the polynomial case, where subcriticality is driven by the dominant contribution of the uniform $l=0$ mode, the subcritical nature here is due to the convex geometry of the $\sinh$. The local asymptotic branches scale as
\begin{equation}
    \psi \sim \pm \sqrt{\frac{\lambda^* - \lambda}{|\lambda_2|}} \left( \sqrt{\frac{105}{16\pi}} \sin\theta \cos^2\theta \cos(2\varphi) \right) + \mathcal{O}(\lambda^* - \lambda).
\end{equation}

\subsection{The Exponential Model}\label{sec:exp}
On a compact manifold like $\mathbb{S}^2$, Gauss's integral constraint \eqref{eq:gauss_constraint} prevents the existence of an exponential model. To ensure a well-posed boundary value problem, we consider the integro-differential Liouville-type equation, coupled with \eqref{eq:gauss_constraint},
\begin{equation}
\begin{cases}
    -\Delta\psi + \lambda e^\psi = \frac{\lambda}{4\pi} \iint_{\mathbb{S}^2} e^\psi \, d\sigma,\\
    \iint_{\mathbb{S}^2} \psi \, d\sigma = 0.
\end{cases}
\end{equation}
Expanding the exponential $$e^\psi = 1 + \psi + \frac{1}{2}\psi^2 + \frac{1}{6}\psi^3 + \mathcal{O}(\psi^4),$$ the $\mathcal{O}(1)$ constants trivially balance across the domain.

\subsubsection*{Order $\mathcal{O}(\varepsilon^1)$} 
At first order, we have,
\begin{equation*}
    -\Delta\psi_1 + \lambda^* \psi_1 = \frac{\lambda^*}{4\pi} \iint_{\mathbb{S}^2} \psi_1\, d\sigma.
\end{equation*}
Imposing the mass constraint yields,
\begin{equation}
    \iint_{\mathbb{S}^2} \psi_1 \, d\sigma = 0.
\end{equation}
Consequently, the non-local linear integral vanishes, reducing the governing equation to 
\begin{equation}
    -\Delta\psi_1 + \lambda^* \psi_1 = 0.
\end{equation}
Because the wave $\psi_1 = Y^*$ belongs to the $l=3$ eigenspace where $\Delta Y^* = -12Y^*$, this drives us to the negative spectrum parameter
\begin{equation}
    \lambda^* = -12.
\end{equation}

\subsubsection*{Order $\mathcal{O}(\varepsilon^2)$}
Let us determine the first-order parameter correction $\lambda_1$. Collecting all $\mathcal{O}(\varepsilon^2)$ terms, the second-order equation reads as
\begin{equation}\label{eq:second_exp}
    L\psi_2 + \lambda_1 Y^* - 6(Y^*)^2 = -\frac{3}{\pi} \iint_{\mathbb{S}^2} \left( \psi_2 + \frac{1}{2}(Y^*)^2 \right) \, d\sigma.
\end{equation}
To isolate $\lambda_1$, we take the inner product with the primary wave $Y^*$. By self-adjointness, the term $\langle L\psi_2, Y^* \rangle$ vanishes.
Furthermore, because the right-hand side is a constant, its inner product with the zero-mean harmonic $Y^*$ ($l=3$) is strictly zero. The projected balance therefore reduces to
\begin{equation}
    \lambda_1 \|Y^*\|_{L^2}^2 - 6 \iint_{\mathbb{S}^2} (Y^*)^3 \, d\sigma = 0.
\end{equation}
Since $Y^*$ is an odd-degree spherical harmonic, its cube $(Y^*)^3$ remains an odd function, forcing its integral over $\mathbb{S}^2$ to vanish. Hence, we obtain
\begin{equation}
    \lambda_1 = 0.
\end{equation}
Recalling that $\iint_{\mathbb{S}^2} \psi_2 \, d\sigma = 0$ by the global mass constraint, the second-order equation becomes
\begin{equation} \label{eq:psi2_equation}
    L\psi_2 = 6(Y^*)^2 - \frac{3}{2\pi} \iint_{\mathbb{S}^2} (Y^*)^2 \, d\sigma.
\end{equation}
Notice that the integral term on the right-hand side is the spatial average (the $l=0$ mode) of $6(Y^*)^2$. Consequently, this term subtracts the mean from the right-hand side, satisfying the mass constraint of $\psi_2$. 
Furthermore, as established in \eqref{eq:c20}, $c_{2,0} = 0$. Therefore, as in \eqref{eq:lambda_2_ris}, inverting the operator $L$ over the orthogonal subspace $l \in \{4, 6\}$ gives
\begin{equation} \label{eq:psi2_solution}
    \psi_2 = -6 \sum_{l \in \{4, 6\}} \sum_{m \in \{0, 4\}} \frac{c_{l,m}}{12 - l(l+1)} R_l^m.
\end{equation}

\subsubsection*{Order $\mathcal{O}(\varepsilon^3)$}
At the third order, the mass constraint gives a constant term $C_3$, such that
\begin{equation}
    L\psi_3 + \lambda_2 \psi_1 - 12\psi_1\psi_2 - 2\psi_1^3 = C_3.
\end{equation}
We project this equation onto the kernel spanned by $\psi_1 = Y^*$. Since $Y^*$ is a zero-mean harmonic ($l=3$), the projection of the right-hand side vanishes ($\langle C_3, Y^* \rangle = 0$). 
Let us emphasize that the mass constraint affects the flow at $\mathcal{O}(\varepsilon^2)$, but not at $\mathcal{O}(\varepsilon^3)$. Exploiting the self-adjointness of $L$ and the normalization $\|Y^*\|_{L^2}^2 = 1$, we isolate $\lambda_2$:
\begin{equation}
    \lambda_2 = 12 \iint_{\mathbb{S}^2} (Y^*)^2 \psi_2 \, d\sigma + 2 \iint_{\mathbb{S}^2} (Y^*)^4 \, d\sigma.
\end{equation}
We now express these integrals in terms of the spectral coefficients $c_{l,m}$ of the harmonic $(Y^*)^2$. By Parseval's identity \eqref{eq:parseval}, the first term on the right-hand side yields the sum of the squared coefficients $\sum c_{l,m}^2$. For the second term, \eqref{eq:psi2_solution} transforms the integral into a spectral series. Then, the bifurcation parameter is given by
\begin{equation}
    \lambda_2 = 12 \sum_{l \in \{4, 6\}} \sum_{m \in \{0, 4\}} \frac{-6 c_{l,m}^2}{12 - l(l+1)} + 2 \sum_{l,m} c_{l,m}^2.
\end{equation}
Notice that,
\begin{equation}
    12 - l(l+1) < 0\quad \text{for all} \quad l \ge 4.
\end{equation}
Thus,
\begin{equation}
    \lambda_2 > 0.
\end{equation}
We conclude that the Liouville bifurcation is \textit{supercritical}, with non-trivial solutions for $\lambda > -12$. The local asymptotic branches scale as
\begin{equation}
    \psi \sim \pm \sqrt{\frac{\lambda - \lambda^*}{\lambda_2}} \left( \sqrt{\frac{105}{16\pi}} \sin\theta \cos^2\theta \cos(2\varphi) \right) + \mathcal{O}(\lambda - \lambda^*).
\end{equation}
Let us remark that the mass constraint acts as a filter that suppresses the $l=0$ zonal mode. By neutralizing the mean-field, the constraint forces the system to rely only on higher-order harmonics, $l\geq4$, shifting the bifurcation from subcritical to supercritical.

\section{Discussion and Conclusions}\label{sec:conclusions}
As summarized in Table \ref{tab:bifurcation_summary}, our spectral analysis suggests the following: while the first-order geometry of tetrahedral non-zonal flows on $\mathbb{S}^2$ is universal, the nature of their topological bifurcation is model-dependent. Bypassing the $(2l+1)$-dimensional kernel degeneracy through equivariant bifurcation theory (Sec.\ref{Sec:Tetra_Group}), we rigorously demonstrated that the bifurcation topology is not an invariant of the spherical Laplacian. Rather, it is determined by the parity and constraints of the nonlinear operator.
\begin{table}[htbp]
\centering
\renewcommand{\arraystretch}{1.6}
\resizebox{\textwidth}{!}{
\begin{tabular}{@{}l l l l l@{}}
\toprule
\textbf{Physical Model} & \textbf{Nonlinearity} $F(\lambda, \psi)$ & \textbf{Nonlinear Interactions} $\mathcal{O}(\varepsilon^2)$ & \textbf{Finite-Amplitude Symmetry} & \textbf{Bifurcation Topology} \\
\midrule
\textbf{Polynomial} & $\mu_1(3\lambda\psi^2 + \psi^3) + \dots$ & $\psi_2 \neq 0$ ($l \in \{0, 4, 6\}$) & Broken ($\psi(-\mathbf{x}) \neq -\psi(\mathbf{x})$) & $\lambda_2 < 0$ (\textbf{Subcritical}) \\
\textbf{Sine-Gordon} & $-\lambda \sin(\psi)$ & $\psi_2 \equiv 0$ (Odd Parity) & Preserved Exactly & $\lambda_2 > 0$ (\textbf{Supercritical}) \\
\textbf{Sinh-Gordon} & $-\lambda \sinh(\psi)$ & $\psi_2 \equiv 0$ (Odd Parity) & Preserved Exactly & $\lambda_2 < 0$ (\textbf{Subcritical}) \\
\textbf{Liouville (Exp)} & $\lambda e^\psi - \frac{\lambda}{4\pi}\iint_{\mathbb{S}^2} e^\psi \, d\sigma$ & $\psi_2 \neq 0$ ($l \in \{4, 6\}$ only) & Broken ($\psi(-\mathbf{x}) \neq -\psi(\mathbf{x})$) & $\lambda_2 > 0$ (\textbf{Supercritical}) \\
\bottomrule
\end{tabular}
}
\vspace{0.2cm}
\caption{Summary of the exact bifurcation structure for the analyzed physical models. The table demonstrates that the bifurcation topology (supercritical or subcritical) is dictated by the $\mathcal{O}(\varepsilon^2)$ nonlinear interactions. Here, $\psi_2$ denotes the second-order spatial correction in the Liapunov-Schmidt asymptotic expansion ($\psi = \varepsilon \psi_1 + \varepsilon^2 \psi_2 + \dots$), representing the secondary flows generated by the primary tetrahedral wave. In particular, the bifurcation direction depends strictly on whether the $l=0$ zonal flow within $\psi_2$ is generated, absent (due to odd parity), or suppressed (due to non-local mass constraints).}
\label{tab:bifurcation_summary}
\end{table}
The most interesting physical consequence of this dependence is the breaking of hemispheric symmetry. The primary instability universally emerges as an odd harmonic ($\psi_1 = Y^*$), enforcing anti-symmetry across the equator. However, for mixed-parity nonlinearities (Polynomial and Liouville models), the quadratic term injects energy into non-trivial nonlinear interactions ($\psi_2 \neq 0$) that rely on even-degree spherical harmonics. This breaks the odd parity ($\psi(-\mathbf{x}) \neq -\psi(\mathbf{x})$), forcing the tetrahedral vortex lattice to distort as the flow gains kinetic energy. Conversely, odd nonlinearities (sine-Gordon and sinh-Gordon) topologically forbid these even-degree nonlinear interactions, preserving the symmetry.

Finally, as shown in Table \ref{tab:bifurcation_summary}, the bifurcation direction is determined by three factors: (i) the cubic term $\psi^3$ in the Taylor expansion of the nonlinearity, (ii) the quadratic interactions $\psi^2$ that drive subcriticality by generating the $l=0$ zonal flow, and (iii) the mass constraint, which suppresses this mean field and forces the system back into a supercritical state.
\newline

\textbf{Conflict of interest.} The author declares that he has no conflict of interest.

\textbf{Data availability.} Data sharing is not applicable.
We do not analyse or generate any datasets, because our work proceeds within a theoretical approach.

\bibliographystyle{plain}
\bibliography{bib}

\end{document}